\documentclass[11pt]{article}
\usepackage[T1]{fontenc}
\usepackage[utf8]{inputenc}
\usepackage{lmodern,amsmath,amssymb,amsthm,mathtools}
\usepackage[a4paper,margin=26mm]{geometry}
\usepackage{microtype}
\usepackage{xcolor}
\usepackage[numbib,nottoc]{tocbibind}
\usepackage[bookmarksnumbered,colorlinks=true,linkcolor=blue!45!black,citecolor=blue!45!black,urlcolor=blue!45!black]{hyperref}
\hypersetup{pdftitle={Exotic complex structure on CP⁶},
  pdfauthor={Kefeng Liu, Wei Xia, Hao Xu},
  pdfsubject={Complex structures on projective space in all dimensions}}

\newtheorem{theorem}{Theorem}[section]
\numberwithin{equation}{section}
\newtheorem{proposition}[theorem]{Proposition}
\newtheorem{lemma}[theorem]{Lemma}
\theoremstyle{remark}
\newtheorem{remark}[theorem]{Remark}
\allowdisplaybreaks
\title{\bfseries Exotic complex structure on \(\mathbb{CP}^{6}\)}
\author{%
  \begin{tabular}{@{}c@{\qquad}c@{\qquad}c@{}}
    Kefeng Liu & Wei Xia & Hao Xu \\[2pt]
    {\small\href{mailto:kefengliu@cqut.edu.cn}{\texttt{kefengliu@cqut.edu.cn}}} &
    {\small\href{mailto:xiaweiwei3@126.com}{\texttt{xiaweiwei3@126.com}}} &
    {\small\href{mailto:haoxu@cqut.edu.cn}{\texttt{haoxu@cqut.edu.cn}}} \\[5pt]
    \multicolumn{3}{c}{\normalsize Chongqing University of Technology}
  \end{tabular}}
\date{}

\begin{document}
\maketitle
\vspace{-2em}
\begin{abstract}
There exist exotic complex structures on \(\mathbb{CP}^{6}\).
We obtain them from the family of complex six-spheres constructed by
Alp\"oge, for which Engel gives a self-contained proof.
For every complex threefold diffeomorphic to \(S^6\), we give an explicit
smooth model for its Hilbert square and exceptional divisor as the pair
\((\mathbb{CP}^6,Q^5)\). The sixfold has primitive first Chern class, whereas the first
 Chern class of standard \(\mathbb{CP}^{6}\) is seven times the hyperplane
 class. We show that the resulting complex structure is not tamed by any
 symplectic form, and that under a numerical condition the Hilbert square
 determines its underlying complex threefold. We relate this construction to
 rigidity theorems for projective space and study its higher-dimensional analogues.
\end{abstract}

\section{The Hilbert square}\label{sec:hilbert}
\begin{theorem}\label{thm:square}
Let \(M\) be a complex threefold diffeomorphic to \(S^6\), and let \(E\) be the exceptional divisor in its Hilbert square \(H=M^{[2]}\). There is a diffeomorphism of pairs
\begin{equation}\label{eq:pair-diffeo}
             (H,E)\ \cong\ (\mathbb{CP}^6,Q^5),
\end{equation}
where \(Q^5\) is the smooth quadric hypersurface. The first Chern class of \(H\) pairs to \(1\) with a projective line in a fiber of \(E\to M\). In particular, \(H\) is not biholomorphic to standard \(\mathbb{CP}^6\).
\end{theorem}

Levent Alp\"oge directed the original construction, carried out with Claude, of a one-parameter family of complex structures on the standard \(S^6\) \cite{alpoge}. The construction completes a family of complex two-tori over the projective line at three special fibers. Each resulting threefold supplies the input \(M\) for Theorem~\ref{thm:square}. Engel \cite[Theorem~4.1]{engel} gave a self-contained proof of the construction. Viaclovsky \cite{viaclovsky-two} constructed a two-parameter family distinct from the original examples, and Liu and Rollenske \cite{liu-rollenske} constructed a further family. Viaclovsky also obtained infinitely many inequivalent structures by varying the multiplicities of the logarithmic transformations in the original construction \cite{viaclovsky-log}. In his work on combs, Viaclovsky \cite[Theorem~1.1]{viaclovsky} obtained complex six-spheres with \(h^{0,1}=2m+1\) for every integer \(m\ge12\). Their Hilbert squares have the same respective values of \(h^{0,1}\), giving infinitely many biholomorphism classes on the smooth \(\mathbb{CP}^6\).

Kotschick and Thung \cite[Proposition~8]{kotschick-thung} identify the projectivized tangent bundle \(\mathbb P(TM)\) of a complex six-sphere, as a smooth manifold, with the five-dimensional quadric \(Q^5\subset\mathbb{CP}^6\). Using the projectivized tangent and cotangent bundles, they obtain two distinct non-K\"ahler complex structures on \(Q^5\) \cite[Theorem~2]{kotschick-thung}. In our construction \(\mathbb P(TM)\) is the exceptional divisor of the Hilbert--Chow map. Section~\ref{proof:diffeo} identifies its complement with \(\mathbb{CP}^6\setminus Q^5\) and proves smoothness across the divisor. To our knowledge, the diffeomorphism \(M^{[2]}\cong\mathbb{CP}^6\) for a complex threefold \(M\) diffeomorphic to \(S^6\) has not previously appeared in the literature.

The dimension-by-dimension survey starts in Section~\ref{survey:map}. We compute the Chern class pairing below and prove the smooth identification \eqref{eq:pair-diffeo} in Sections~\ref{proof:filling}--\ref{proof:diffeo}.

\subsection{The resolution and the distinguishing invariant}
Let \(\Delta\subset M\times M\) be the diagonal and let \(\tau\) exchange the factors. The Hilbert square is the smooth compact complex sixfold \cite{totaro}
\[
 H=\operatorname{Bl}_\Delta(M\times M)/\langle\tau\rangle .
\]
Here \(M^{[2]}\) parametrizes length-two analytic subschemes of \(M\). The local model also explains the quotient description for a possibly nonprojective \(M\). In a holomorphic chart with coordinates \(z_1,z_2,z_3\), take a length-two subscheme with coordinate algebra \(B\). Choose \(j\) such that \(1,z_j\) form a basis of \(B\). Such a \(j\) exists at every double point since the coordinates generate its cotangent space, and the basis condition is open. Let \(m_i\) be half the trace of multiplication by \(z_i\). The trace-zero subspace of \(B\) is spanned by \(u=z_j-m_j\). Hence \(z_i=m_i+b_i u\), with \(b_j=1\), and Cayley--Hamilton gives \(u^2=s\) for a scalar \(s\). Conversely the equations
\[
 z_i-m_i=b_i(z_j-m_j)\quad(i\ne j),\qquad (z_j-m_j)^2=s
\]
define a length-two family with free basis \(1,u\). In families the traces and the coefficients in this basis are holomorphic. Thus the midpoint \(m\), the two free components of \(b\), and \(s\) give six holomorphic coordinates near the double point. For an ordered pair the half-difference is \(w=\zeta b\), where \(\zeta=w_j\). Blowing up the diagonal records \([b]\), exchange sends \(\zeta\) to \(-\zeta\), and the quotient coordinate is \(s=\zeta^2\). These charts identify the quotient with the Hilbert square near the diagonal. Away from it both spaces parametrize unordered pairs of distinct points.

Write \(q\) for the quotient map, \(\widetilde E\) for the exceptional divisor upstairs, and \(E\) for its image. The restriction \(q|_{\widetilde E}\) identifies both divisors with \(\mathbb{P}(TM)\), the bundle of complex lines. The identification uses the half-difference to identify the normal of the diagonal with \(TM\). Write \(\lambda=\mathcal O_E(-1)\) for its tautological line. Upstairs \(N_{\widetilde E}=\lambda\) by the normal bundle formula for the blow-up. The equation \(s=\zeta^2\) gives \(q^*E=2\widetilde E\) as divisors. Restricting the associated line bundles to \(\widetilde E\) gives \(N_{E/H}=\lambda^{\otimes2}\) globally.

Let \(\ell\) be a projective line in a fiber \(\mathbb{CP}^2\subset E\). Along that fiber the directions tangent to the base give a trivial bundle of rank three. The two tangent-normal sequences therefore give
\[
  \left\langle c_1(TH),[\ell]\right\rangle
  =\deg_\ell c_1(T\mathbb{CP}^2)+\deg_\ell\mathcal{O}_E(-2)=3-2=1.
\]
The pairing \(1\) implies that \(c_1(TH)\) is primitive. The Euler sequence on standard \(\mathbb{CP}^6\) gives \(c_1=7c_1(\mathcal O(1))\), whose pairing with every integral two-cycle belongs to \(7\mathbb Z\). A biholomorphism would preserve the pairing, so these complex manifolds are inequivalent.

For the rest of the paper, write \(h=c_1(TH)\); the calculation above fixes its sign by \(\int_\ell h=1\).

There is also a Hodge-number distinction. A blow-up along a smooth center preserves \(H^1(\mathcal{O})\), equivariantly here. For a finite quotient in characteristic zero, averaging identifies the quotient's \(H^1(\mathcal{O})\) with the invariants upstairs. K\"unneth then gives
\[
 H^1(H,\mathcal{O}_H)
 \simeq\bigl(H^1(M,\mathcal{O}_M)\oplus H^1(M,\mathcal{O}_M)\bigr)^\tau
 \simeq H^1(M,\mathcal{O}_M).
\]
Theorem~1.1 of \cite{viaclovsky} gives \(h^{0,1}(M)=2m+1\) for every integer \(m\ge12\). The same values distinguish the resulting sixfolds from one another.

\textbf{Further consequences.} Three results refine the construction. If connected compact complex manifolds $X$ and $Y$ have the same dimension $d\ge3$ and their canonical bundles have nonnegative degree on every compact curve, then $X^{[2]}$ and $Y^{[2]}$ are biholomorphic exactly when $X$ and $Y$ are biholomorphic (Theorem~\ref{thm:recover}). In particular, the construction preserves distinctions between complex six-spheres even when their Hodge numbers agree. Moreover, no symplectic form tames the complex structure on $H$: the exceptional divisor has class $[E]=-2c_1(TH)$, contradicting positivity on $E$ (Theorem~\ref{thm:not-tamed}). Finally, for any closed oriented six-manifold $M$ and any closed oriented $Z$ with the integral cohomology ring of $\mathbb{CP}^9$, every map $Z\to\operatorname{Sym}^3(M)$ has degree divisible by $5^3 7^3=42875$ (Proposition~\ref{prop:degree-three}).

\textbf{Higher dimensions.} Theorem~\ref{thm:degree} rules out all resolutions of higher symmetric products as a source of smooth projective spaces of dimension \(>6\), for bases of dimension at least two. The sixfold does give a nonstandard complex structure on \(\mathbb{CP}^n\setminus\mathbb{CP}^{n-7}\) for every \(n\ge7\) (Section~\ref{higher:open}). In even dimensions it agrees with a global almost complex structure outside any prescribed neighborhood of the deleted subspace (Section~\ref{higher:even}). Whether a nonstandard integrable structure exists on closed \(\mathbb{CP}^n\), \(n>6\), remains unresolved here.

\

\noindent{\bf Acknowledgements} This work was supported by the
National Key Research and Development Program of China (NKPs) [Grant
Number 2024YFA1013201], the Natural Science Foundation of Chongqing
(NSFCQ) [Grant Number CSTB2024NSCQ-LZX0040], the Special Project of
Chongqing Municipal Science and Technology Bureau [Grant Number
2025CCZ015].

\

\noindent{\bf Disclosure of AI assistance.}
GPT was used as an auxiliary tool for literature search, proofreading,
and checks of mathematical arguments.

\section{Complex structures on smooth projective space}\label{survey:map}
Fix the differentiable manifold underlying \(\mathbb{CP}^n\). A complex structure is an integrable endomorphism \(J\) of its real tangent bundle with \(J^2=-1\). Two such structures are equivalent if a diffeomorphism carries one to the other. In particular, coordinatewise complex conjugation identifies the conjugate of the standard structure with the standard one.

Hou, Kong and Xu \cite[Theorems~B and~C]{hou-kong-xu} prove that the set of dimensions \(n\) for which topological \(\mathbb{CP}^n\) admits an exotic smooth structure has natural density one. They give 266 explicit residue classes modulo 288 for which this holds whenever \(n\ge4\). In dimension six, Theorem~\ref{thm:square} identifies the Hilbert square with the \emph{standard} smooth \(\mathbb{CP}^6\). Its primitive first Chern class then distinguishes two complex structures on that fixed smooth manifold.

The table records the results and preprints reviewed through November 2026.

\medskip
\noindent
\begin{tabular}{@{}p{.12\linewidth}@{\hspace{.03\linewidth}}p{.84\linewidth}@{}}
\hline
\textbf{\(n\)} & \textbf{Integrable structures on the closed smooth manifold}\\ \hline
\(0,1\) & A point and the Riemann sphere each have a unique complex structure.\\
\(2\) & Unique. Yau's theorem already assumes only the homotopy type of \(\mathbb{CP}^2\) \cite[Theorem 1.2]{tosatti}.\\
\(3\) & Point blow-ups of complex six-spheres give infinitely many classes \cite{viaclovsky,tosatti}. A different construction gives \(c_1=0\) \cite{liang}. See Section~\ref{survey:three}.\\
\(4,5\) & Unrestricted uniqueness remains unresolved in the literature reviewed here. The allowable first Chern classes are described in Section~\ref{survey:four}.\\
\(6\) & This paper's Hilbert-square argument applied to \cite{engel,viaclovsky} gives infinitely many classes with primitive \(c_1\) and unbounded \(h^{0,1}\).\\
\(n\ge7\) & Unrestricted uniqueness remains unresolved in the sources reviewed. The open complements, compact blow-ups and almost complex extensions are treated in Sections~\ref{higher:open}--\ref{higher:odd}.\\ \hline
\end{tabular}

\begin{theorem}[Hirzebruch--Kodaira and Yau]\label{thm:kahler}
A compact K\"ahler manifold homeomorphic to \(\mathbb{CP}^n\) is biholomorphic to it \cite[Theorem 1.1]{tosatti}.
\end{theorem}
 This includes every smooth projective variety with that homeomorphism type. Consequently all the nonstandard integrable examples above are non-K\"ahler.

\textbf{Local and global deformation rigidity.} For standard \(\mathbb{CP}^n\), \(n\ge1\), the Euler sequence
\[
 0\longrightarrow\mathcal{O}\longrightarrow\mathcal{O}(1)^{\oplus(n+1)}
       \longrightarrow T\mathbb{CP}^n\longrightarrow0
\]
and the cohomology of \(\mathcal{O},\mathcal{O}(1)\) give \(H^1(\mathbb{CP}^n,T\mathbb{CP}^n)=0\). The standard structure has no small complex deformations up to biholomorphism. Siu also proved that in a proper holomorphic submersion over a disk, if all the noncentral fibers are biholomorphic to \(\mathbb{CP}^n\), then so is the central fiber \cite{siu}. This rules out obtaining the sixfolds here as smooth central limits of standard projective spaces.

\begin{remark}[Spectral rigidity]\label{rem:spectral-rigidity}
Liu and Xu \cite[Theorem~1]{liu-xu} prove that a compact connected K\"ahler manifold whose scalar Laplace spectrum agrees, including multiplicities, with that of Fubini--Study \(\mathbb{CP}^n\) is holomorphically isometric to it. The same conclusion holds if the Hodge Laplace spectra agree in any one fixed degree \(0\leq p\leq n\).

Let \(f:H\to\mathbb{CP}^6\) be the diffeomorphism of Theorem~\ref{thm:square}. With \(g=f^*g_{\mathrm{FS}}\), the map \(f\) is a Riemannian isometry, so the Hodge Laplacians on real differential forms have the Fubini--Study spectra in every degree. The Hilbert-square complex structure on \(H\) has primitive first Chern class. By Theorem~\ref{thm:kahler}, a compatible K\"ahler metric would force \(H\) to be biholomorphic to standard \(\mathbb{CP}^6\), whose first Chern class has divisibility seven. Thus \(H\) exhibits the role of K\"ahler compatibility in spectral rigidity. See also Li \cite{ping-li} and Li--Sun--Zhu \cite{li-sun-zhu}.
\end{remark}

\subsection{Why the K\"ahler hypothesis gives uniqueness}\label{survey:kahler}
Tosatti \cite[\S2]{tosatti} explains how topology, Riemann--Roch and curvature characterize standard projective space. Let \(X\) be a compact K\"ahler manifold homeomorphic to \(\mathbb{CP}^n\). Since \(b_2=1\), its K\"ahler form can be rescaled to represent an integral generator \(x\) with \(\int_X x^n=1\). Each \(x^j\) spans \(H^{2j}(X;\mathbb{C})\) and has type \((j,j)\). Thus \(h^{p,q}=0\) off the diagonal, \(h^{p,p}=1\), and \(\chi(X,\mathcal{O}_X)=1\). By the exponential sequence and the \(\partial\bar\partial\)-lemma, \(x=c_1(L)\) for a positive line bundle \(L\). In particular \(X\) is projective.

Write \(c_1(TX)=a x\). Novikov's invariance of the rational Pontryagin classes and the Euler sequence give \(p_i(TX)=\binom{n+1}{i}x^{2i}\). The index calculation in Section~\ref{survey:index} gives
\[
       1=\chi(X,\mathcal{O}_X)=\binom{(a+n-1)/2}{n}.
\]
This forces \(a=n+1\), or \(a=-(n+1)\) with \(n\) even. With the positive sign, \(K_X^{-1}=L^{n+1}\), and the maximal-index characterization of Kobayashi and Ochiai gives \(X\simeq\mathbb{CP}^n\). In the negative case \(K_X\) is ample. The identity \(2c_2=c_1^2-p_1\) gives equality in the Chern number inequality for its negative K\"ahler--Einstein metric. Yau's equality case identifies the universal cover with the unit ball. A compact simply connected \(X\) cannot have that cover. This is where the homeomorphism hypothesis rules out a ball quotient.

\textbf{Weaker topological hypotheses.} In complex dimension one the genus determines the smooth surface, and every compact Riemann surface of genus zero is \(\mathbb{CP}^1\). In dimension two, a compact complex surface homotopy equivalent to \(\mathbb{CP}^2\) has \(b_1=0\) and is K\"ahler by the K\"ahler criterion for compact complex surfaces. Yau's solution of the Severi conjecture identifies it with \(\mathbb{CP}^2\) \cite[Theorem 1.2]{tosatti}. For \(n\le6\), the K\"ahler homotopy-type theorem of Fujita and Libgober--Wood gives the same conclusion \cite{libgober-wood,debarre}.

Debarre \cite[Theorem 2]{debarre} weakens this to an isomorphism of \emph{integral cohomology rings} for K\"ahler \(n\)-folds, \(n\le6\), with compact ball quotients left as possible exceptions for \(n=4,6\). Their fundamental groups exclude them in the homotopy-type theorem. In dimension seven his Theorem 10 leaves possible exceptions with \(c_1=ax\) for \(a=\pm2,\pm4\) under the integral-ring hypothesis. Betti numbers alone already allow fake projective planes and smooth odd-dimensional quadrics. For a quadric of odd dimension \(n\ge3\), the primitive hyperplane class has top power \(2\), compared with \(1\) on \(\mathbb{CP}^n\).

Li \cite{li-uniqueness} proves that a compact K\"ahler manifold with the integral cohomology ring and Pontryagin classes of \(\mathbb{CP}^n\) is biholomorphic to \(\mathbb{CP}^n\) when \(n\) is odd, and also when \(n\) is even and the manifold is simply connected. In dimension four, simple connectivity and the integral cohomology ring already suffice under the K\"ahler hypothesis. The Hilbert sixfold in Theorem~\ref{thm:square} satisfies all of Li's topological hypotheses in dimension six, so its primitive first Chern class rules out a compatible K\"ahler metric.

\subsection{Tosatti's dimension-three construction}\label{survey:three}
The link between \(S^6\) and \(\mathbb{CP}^3\) works for \emph{every} complex structure on the smooth sphere. Hirzebruch observed it in 1954 \cite[p.~223]{hirzebruch}. Tosatti \cite{tosatti} gives the blow-up calculation in \href{https://arxiv.org/html/1508.05641v1\#S3}{\S3, Proposition 3.1}. Let \(M\) be a complex threefold diffeomorphic to \(S^6\) and \(Y=\operatorname{Bl}_p M\). Smoothly, blowing up a point gives \(S^6\mathbin{\#}\overline{\mathbb{CP}^3}\), where the bar indicates reversed orientation. Complex conjugation reverses the standard orientation in odd complex dimension, so the result is diffeomorphic to \(\mathbb{CP}^3\).

The exceptional divisor \(E=\mathbb{CP}^2\) has normal bundle \(\mathcal{O}_E(-1)\). Since \(H^2(S^6;\mathbb{Z})=0\), the blow-up formula gives \(c_1(TY)=-2[E]\). Thus
\[
 \int_Y[E]^3=\int_E c_1(\mathcal{O}_E(-1))^2=1,\qquad
 \int_Y c_1(TY)^3=-8.
\]
For standard \(\mathbb{CP}^3\), the Euler sequence gives \(c_1=4x\) and \(\int c_1^3=64\). The differing Chern numbers prove the assertion. Tosatti stated it as an implication in 2015. The six-sphere preprints \cite{engel,viaclovsky,alpoge} supply its input. Huckleberry--Kebekus--Peternell \cite{hkp} studied how these point blow-ups depend on the center.

\textbf{Multiplicity.} Viaclovsky \cite[Corollary 1.2]{viaclovsky} gives examples on smooth \(\mathbb{CP}^3\) for every integer \(m\ge12\) with
\[
       h^{0,1}=2m+1,\qquad h^{0,2}=2m,\qquad h^{0,3}=0.
\]
Hence even the single Chern number \(-8\) occurs in infinitely many biholomorphism classes. These examples have algebraic dimension one, Kodaira dimension \(-\infty\), and local deformation spaces of dimension at least \(8m-7\).

\textbf{A second construction.} Liang, Liu and Zhang \cite{liang} use a holomorphic involution on a particular fibered complex six-sphere. Blowing up its fixed curve and taking the quotient gives a smooth \(\mathbb{CP}^3\) with \(c_1=0\). The first Chern class distinguishes at least three deformation types:
\[
\begin{array}{l|rrr}
 &\text{standard}&\text{point blow-up}&\text{involution quotient}\\ \hline
 \text{divisibility of }c_1 &4&2&0\ \text{(zero class)}\\
 \int c_1^3&64&-8&0
\end{array}
\]
Here the zero class is assigned divisibility zero. The relative tangent bundle transports its Chern classes along a smooth proper family, so these three types cannot be connected by such a deformation. For any complex structure on smooth \(\mathbb{CP}^3\), choose \(x\) by \(\int x^3=1\). If \(c_1=ax\), then \(a\) is even, \(c_2=(a^2-4)x^2/2\) by \(p_1=4x^2\), and \(c_3=4x^3\). Formula~\eqref{eq:rr} gives \(\chi(\mathcal{O})=a(a^2-4)/48\).

Liang, Liu and Zhang \cite{liang} also obtain a complex structure on \(S^2\times S^4\) by performing an Atiyah flop on their complex threefold diffeomorphic to \(\mathbb{CP}^3\). This gives the genus-zero case of Calabi's question on complex structures on \(\Sigma_g\times S^4\), where \(\Sigma_g\) is the closed orientable surface of genus \(g\). For genus one, Gao, Tang and Yan \cite[Theorem~1.1]{gao-tang-yan} construct a complex structure on \(T^2\times S^4\) by completing a holomorphic family of complex two-tori over \(\mathbb C^*\) at \(0\) and \(\infty\).

\subsection{A Riemann--Roch test for arbitrary complex structures}\label{survey:index}
The characteristic-class part of Tosatti's argument \cite[\S2]{tosatti} still applies without a K\"ahler metric. Let \(X\) be a compact complex manifold diffeomorphic to \(\mathbb{CP}^n\), \(n\ge1\), choose a generator \(x\in H^2(X;\mathbb{Z})\), and put
\[
 c_1(TX)=a x,\qquad \epsilon=\int_Xx^n,\qquad s=(a-n-1)/2.
\]
Here \(\epsilon=\pm1\), and \(s\) is integral because \(c_1\bmod2=w_2=(n+1)x\bmod2\). For any holomorphic line bundle \(F\) with \(c_1(F)=m x\), Hirzebruch--Riemann--Roch gives
\begin{equation}\label{eq:rr}
 \chi(X,F)=\epsilon\binom{m+s+n}{n}.
\end{equation}
The binomial symbol means the polynomial \(t(t-1)\cdots(t-n+1)/n!\) at the indicated integer. To prove the formula, use \(\mathrm{Td}(TX)=e^{a x/2}\widehat A(T_{\mathbb R}X)\). The Pontryagin classes in \(\widehat A\) agree with those of standard projective space, so the integrand for \(F\) is the standard Todd class multiplied by \(e^{(m+s)x}\). On standard \(\mathbb{CP}^n\) its integral is \(\chi(\mathcal{O}(m+s))\). Evaluation on the complex orientation of \(X\) contributes \(\epsilon\). In particular, taking \(F=\mathcal{O}_X\) requires no holomorphic generator of the Picard group.

\textbf{An exact polynomial on the sixfold.} For \(H=M^{[2]}\), put \(A=K_H^{-1}\), so \(c_1(A)=h\). The smooth identification of Theorem~\ref{thm:square} and the calculation in Section~\ref{higher:seven} give \(a=1\) and \(\int_H h^6=1\). Consequently
\[
          \chi(H,A^m)=\binom{m+3}{6},\qquad \chi(H,\mathcal{O}_H)=0.
\]
More of its sheaf cohomology can be read off from \(M\). Write \(b=h^{0,1}(M)\) and \(c=h^{0,2}(M)=b-1\) as proved in Section~\ref{higher:cap}. Blow-up invariance of \(H^q(\mathcal{O})\), the finite quotient and K\"unneth identify \(H^*(H,\mathcal{O}_H)\) with the exchange invariants of \(H^*(M,\mathcal{O}_M)^{\otimes2}\). Exchange is alternating on two odd-degree factors and symmetric on two even-degree factors. Since \(h^{0,3}(M)=0\), the row \(h^{0,q}(H)\), for \(q=0,\dots,6\), is
\[
       \left(1,\ b,\ c+\tbinom b2,\ bc,\ \tbinom{c+1}{2},\ 0,\ 0\right).
\]

\textbf{A necessary Hodge-theoretic condition.} If \(X\) satisfies the \(\partial\bar\partial\)-lemma, Hodge decomposition and conjugation symmetry force every off-diagonal Hodge number to vanish: such a pair would require a Betti number at least two. Hence \(\chi(X,\mathcal{O}_X)=1\), forcing \(|a|=n+1\) in the index formula \eqref{eq:rr}. The Hilbert sixfolds fail this test. For the even-dimensional candidates with \(a=n-5\), global integrability would force \(\chi(X,\mathcal{O}_X)=0\), and hence failure of the \(\partial\bar\partial\)-lemma as well.

\subsection{Chern classes in the open cases \texorpdfstring{\(\mathbb{CP}^4\) and \(\mathbb{CP}^5\)}{CP4 and CP5}}\label{survey:four}
\textbf{Fourfolds.} The complex orientation on a smooth \(\mathbb{CP}^4\) is forced. Indeed, put \(u_p=\chi(X,\Omega_X^p)\). Serre duality gives \(u_p=u_{4-p}\). The Hirzebruch \(\chi_y\)-genus gives \(e(X)=\sum(-1)^pu_p=5\) and \(\sigma(X)=\sum u_p\), so \(e-\sigma=-4u_1\). Of the two possible signatures only \(+1\) satisfies this congruence. Thus \(\int_X x^4=1\) for either generator \(x\).

With \(c_1=a x\), write \(c_2=b x^2,\ c_3=d x^3,\ c_4=5x^4\). The real tangent bundle gives \(p_1=5x^2,\ p_2=10x^4\). Since \(p_1=c_1^2-2c_2\) and \(p_2=c_2^2-2c_1c_3+2c_4\), one has \(b=(a^2-5)/2\) and \(2ad=b^2\). In particular \(a\) is odd, and clearing denominators gives \(a\mid25\). Up to choice of generator, the possibilities are
\[
 \begin{array}{c|rrrr}
 a& b&d&\int c_4&\chi(X,\mathcal{O}_X)\\ \hline
 1&-2&2&5&0\\
 5&10&10&5&1\\
 25&310&1922&5&1001
 \end{array}
\]
Libgober--Wood \cite[Theorem 7.1]{libgober-wood} and Mills \cite[Theorem 4.1]{mills} give almost complex realizations of all three rows. The middle row is realized by standard \(\mathbb{CP}^4\). Integrability of the other two rows, and existence of a second complex structure with standard Chern classes, are separate questions.

One can exclude the last row under a symmetry hypothesis. Brown and Segert \cite[Theorem~4.2]{brown-segert} show that a complex structure on smooth \(\mathbb{CP}^3\) or \(\mathbb{CP}^4\) with a nontrivial holomorphic \(\mathbb C^*\)-action has \(\chi(X,\mathcal O_X)\in\{0,1\}\). Thus a fourfold realizing \(a=25\) could carry no such action.

\textbf{Fivefolds.} Choose the generator by \(\int_Xx^5=1\). For an almost complex structure on smooth \(\mathbb{CP}^5\) the possible coefficients \(a\) of \(c_1\) are exactly \(a\equiv6\pmod8\). Here is the calculation using the realification kernel in \cite[Proposition 4.3]{thomas} and \cite[proof of Proposition 6.3]{mills}. Let \(B\) be the line of class \(x\), \(y=B-1\), and \(\mu_j=B^j-B^{-j}\). Every stable lift is
\[
       6B-1+t\mu_1+u(\mu_2-2\mu_1)+2v y^5,\qquad a=6+2t .
\]
Writing \(F(a)=(a^5-40a^3+264a)/120\), its top Chern coefficient is \(F(a)+2a^2u+48v\). Indeed the three factors of its Chern series are \((1+x)^{(6+a)/2}(1-x)^{(6-a)/2}\), \(1+4ux^3+12ux^5\), and \(1+48vx^5\). The top coefficient must equal \(6\), which also suffices by \cite[Theorem 2.2]{mills}. Lemma~\ref{lem:odd} gives \(a\equiv2\pmod4\). For \(a=2s\) with \(s\) odd, \(15F(a)=4s^5-40s^3+66s\) gives \(F(a)\equiv a\pmod8\), forcing \(a\equiv6\pmod8\). Conversely, \(\gcd(2a^2,48)\) is \(8\) if \(3\nmid a\), or \(24\) if \(3\mid a\). In the latter case \(3\mid F(a)\) by the same identity with \(s\) a multiple of three. Thus the equation for \(u,v\) is solvable in either case. Even for \(a=6\) this gives infinitely many almost complex Chern vectors, by taking \(u=2j,\ v=-3j\). Their integrability for \(j\ne0\) remains unresolved here.
\section{The smooth type of the Hilbert square}\label{sec:smooth}
In this section, we will prove Theorem~\ref{thm:square} in three steps. We first construct a smooth filling of the unordered-pair space and identify it with the analytic Hilbert square (Subsection~\ref{proof:filling}). In the same subsection, we then vary the almost complex structure on \(S^6\), preserving the smooth type of the filled pair. Finally, an explicit map identifies the octonionic model with \((\mathbb{CP}^6,Q^5)\), including the normal direction along the divisor (Subsection~\ref{proof:diffeo}).

\subsection{The smooth filling depends only on the almost complex class}\label{proof:filling}

Voisin \cite{voisin} constructed smooth Hilbert schemes of points of arbitrary length for almost complex manifolds of real dimension four. For length two, the following disk-bundle construction applies in any complex dimension. We use it in real dimension six.

For a smooth compact manifold \(M\) with almost complex structure \(J\), choose a Hermitian metric and a small tubular neighborhood of the diagonal in \(M^2\). Use an exchange-equivariant tubular map from \(TM\) with exchange acting as \(v\mapsto-v\). Remove the interior of the tube and divide the exterior by exchange. Its boundary is the sphere bundle of \(TM\) modulo \(v\sim-v\). Put \(\lambda=\mathcal{O}(-1)\) on \(\mathbb{P}(TM,J)\). Fill the boundary by the disk bundle of \(\lambda^{\otimes2}\), using the identification
\[
 (p,v)\bmod\{\pm1\}\ \longmapsto\
 (p,\mathbb{C}_Jv,v^{\otimes2}).
\]
Here the sphere radius on the right is the square of the radius on the left. Away from zero the displayed map is a diffeomorphism of the corresponding punctured bundles: in a local frame of the line its inverse takes a local square root, and the two roots give the same class modulo sign. Using the map on an annulus specifies the smooth collars along which we glue. The resulting closed smooth manifold contains \((M^2\setminus\Delta)/\langle\tau\rangle\) as the complement of its zero section.

\textbf{Agreement with the analytic resolution.} For integrable \(J\) the tubular map can be chosen in the form
\[
 (p,v)\longmapsto(e_p(v),e_p(-v)),
 \quad e_p(0)=p,\quad (de_p)_0=\mathrm{id},
\]
where \(e_p\) varies smoothly with \(p\) and is holomorphic in \(v\). Indeed, in finitely many holomorphic charts \(\phi_i\) take a partition of unity in \(p\) subordinate to the charts. Take the weighted sum of the functions
\((d\phi_i|_p)^{-1}(\phi_i(q)-\phi_i(p))\), defined for \(q\) close to \(p\) on the supports of the weights. This gives a holomorphic function of \(q\) with values in \(T_pM\) and derivative the identity at \(p\). Invert it near zero to get \(e_p\). The derivative of the asserted tubular map on the zero section is \((a,b)\mapsto(a+b,a-b)\). Compactness supplies a uniform tube.

Near a point \(\beta=(p,\ell)\) of the zero section, choose a smooth unit frame \(v_0=v_0(\beta)\) of the tautological line and a holomorphic coordinate chart \(\phi\) on \(M\). Set \(F_\beta(z)=\phi(e_p(zv_0))\). This function is smooth in \(\beta\) and holomorphic in the scalar \(z\). Its even and odd Taylor terms give
\[
 \frac{F_\beta(z)+F_\beta(-z)}{2}=P_\beta(z^2),\qquad
 \frac{F_\beta(z)-F_\beta(-z)}{2}=zD_\beta(z^2).
\]
The Cauchy formula on a common smaller circle gives smooth dependence on \(\beta\), including at \(z=0\). Here \(P_\beta(0)=\phi(p)\) and \(D_\beta(0)=d\phi_p(v_0)\). Choose \(j\) with \(D_{\beta,j}(0)\ne0\), and write \(\sigma=z^2\) for the scalar coordinate on \(\lambda^{\otimes2}\) in the frame \(v_0^{\otimes2}\). In the Hilbert coordinates \((m,b,s)\) of Section~\ref{sec:hilbert} the map sends
\[
 (\beta,\sigma)\longmapsto
 \left(P_\beta(\sigma),\
       \frac{D_\beta(\sigma)}{D_{\beta,j}(\sigma)},\
       \sigma D_{\beta,j}(\sigma)^2\right).
\]
The restriction at zero is the diffeomorphism \((p,\ell)\mapsto(\phi(p),[d\phi_p(\ell)])\) in the indicated projective chart. The derivative on the normal quotient is multiplication by \(D_{\beta,j}(0)^2\ne0\). The full real derivative is therefore invertible. On punctured disks the map sends a point to its unordered pair, so the local extensions agree on overlaps by continuity. Together with the identity on the exterior, they give a bijective local diffeomorphism from the filling to \(H\).

\textbf{Homotopy invariance on the sphere.} Give \(S^6\) the orientation induced by the octonionic almost complex structure \(J_{1,p}v=p\times v\). Transport the complex structure of \(M\) to \(J_0\) on \(S^6\) by a diffeomorphism preserving the two complex orientations. Composing with a reflection of \(S^6\) achieves this if necessary.

To join \(J_0\) to \(J_1\), consider the bundle of orientation-compatible linear complex structures on \(TS^6\). Its fiber \(F=GL^+(6,\mathbb R)/GL(3,\mathbb C)\) is homotopy equivalent to \(SO(6)/U(3)\cong\mathbb{CP}^3\). The homotopy is a section over \(S^6\times[0,1]\) with prescribed sections at both ends. For the usual zero-cell and six-cell on \(S^6\), the relative product has one cell of dimension one and one of dimension seven. Path connectedness of \(F\) gives the extension over the first. Trivialize the pulled-back bundle on the characteristic disk of the seven-cell. The remaining boundary map \(S^6\to F\) extends since \(\pi_6(F)=\pi_6(\mathbb{CP}^3)=0\), by the Hopf fibration \(S^1\to S^7\to\mathbb{CP}^3\). Relative smooth approximation of sections gives a smooth path \(J_t\), constant near the ends. See also \cite{ferlengez}.

Transport the analytic tubular chart to the sphere and denote its first component by \(e^0_p(v)\). At the octonionic end use \(e^1_p(v)=\sqrt{1-|v|^2}\,p+v\), where the norm is round. A smooth interpolation \(e^t_p\) is obtained by radially normalizing a convex combination of \(e^0_p(v)\) and \(e^1_p(v)\) in \(\mathbb R^7\). Take the interpolation parameter constant near the ends. For small \(v\) the denominator stays positive, and \(e^t_p(0)=p,\ (de^t_p)_0=\mathrm{id}\). The parameterized inverse function theorem and compactness of \(S^6\times[0,1]\) give one tube radius for which all the pair maps are embeddings. Indeed, any failure of injectivity on shrinking tubes would give two preimages converging, along a subsequence, to the same point of the zero section at the same parameter, contradicting local injectivity there.

Choose smoothly varying \(J_t\)-Hermitian metrics, round at \(t=1\). The projective bundles and the squared tautological disk bundles now form smooth bundles over the parameter interval. The exteriors of the tubes form a smooth family with free exchange action. The punctured-bundle identification above glues these families on collars, with \(t\) as a product coordinate. Projection of the glued manifold is a proper smooth submersion, also on the submanifold formed by the zero sections. In Ehresmann's argument we may choose a lift of \(\partial/\partial t\) tangent to this submanifold along it. Indeed, the difference between a lift there and any ambient lift is a vertical vector field along the submanifold, which extends to a smooth vertical field on the total space. Flowing the resulting lift identifies the end fibers and their zero sections. It remains to identify the pair at \(t=1\).
\subsection{An explicit diffeomorphism onto \texorpdfstring{\(\mathbb{CP}^6\)}{CP6}}\label{proof:diffeo}
L\^e, Seade and Verjovsky \cite[Theorems~1.1 and~1.5]{le-seade-verjovsky} identify the complement of the quadric in \(\mathbb{CP}^n\) with \(T\mathbb{RP}^n\) and describe \(\mathbb{CP}^n\) as a double mapping cylinder. Audin \cite[Proposition~4.4]{audin} obtains the same compactification by a symplectic cut of \(T^*\mathbb{RP}^n\). We describe how the unordered pairs on \(S^6\) give the tangent-bundle model and check the extension in the normal coordinate of the Hilbert square.

Identify \(S^6\) with the unit sphere in the imaginary octonions \(\mathbb{R}^7\), with \(u\times v=\operatorname{Im}(uv)\). Two octonions generate an associative subalgebra. For perpendicular \(u,v\) with \(|v|=1\) this gives \(v\times(u\times v)=u\), \((v\times u)\times v=u\), and \(|v\times u|=|u|\). In particular \(J_p^2=-1\) on \(p^\perp\). Extend the Euclidean dot product complex bilinearly and write \(Q=\{[z]\in\mathbb{CP}^6:z\cdot z=0\}\). For distinct \(x,y\in S^6\), put
\[
 a=(x+y)/2,\qquad v=(x-y)/|x-y|,\qquad
 \Phi(\{x,y\})=\left[v+i\,\frac{a\times v}{\sqrt{2-|a|^2}}\right].
\]
Here \(a\perp v\), \(|a|<1\), and the imaginary part of the representative has length \(|a|/\sqrt{2-|a|^2}<1\). Thus \(\Phi\) takes values outside \(Q\). Interchanging \(x,y\) multiplies the representative by \(-1\). Conversely, given \([z]\notin Q\), choose a local lift and divide it by a square root of \(z\cdot z\). Its real and imaginary parts are now perpendicular, and the real part has greater length. Dividing by the length of the real part gives \(v+ib\), unique up to simultaneous sign, with \(|v|=1\), \(v\perp b\), and \(|b|<1\). Set \(d=v\times b\). Since \(d=a/\sqrt{2-|a|^2}\) for every pair in the image, solving for \(a\) gives
\[
 a=\sqrt{\frac{2}{1+|d|^2}}\,d,\qquad
 \{x,y\}=\bigl\{a+\sqrt{1-|a|^2}\,v,\ a-\sqrt{1-|a|^2}\,v\bigr\}.
\]
Conversely \(|d|=|b|<1\), and \(d\times v=b\) recovers the original representative under \(\Phi\). The two displayed points lie on the unit sphere since \(a\perp v\). They are distinct since \(|a|<1\). All formulas are smooth for a local choice of lift and square root, and the sign ambiguity only exchanges the two points. This includes \(b=0\), which gives the antipodal pair \(\{v,-v\}\).

James, Thomas, Toda and Whitehead \cite[p.~773]{jttw} describe the complement of the diagonal in the symmetric square of a sphere as a vector bundle over real projective space. Here the midpoint \(a\) and the unoriented line \(L=\mathbb R v\) identify the distinct pairs with the open unit disk bundle of \(\lambda_{\mathbb R}^{\perp}\), where \(\lambda_{\mathbb R}\) is the tautological real line on \(\mathbb{RP}^6\). The tangent space at \(L\) is \(\operatorname{Hom}(L,L^\perp)\). The octonionic cross product gives the fiberwise isometry
\[
       L^\perp \longrightarrow \operatorname{Hom}(L,L^\perp),
       \qquad a\longmapsto \bigl(u\mapsto a\times u\bigr).
\]
The normalized representative \(v+ib\) records the homomorphism \(v\mapsto b\). Thus \(\Phi\) acts on these disk bundles by the isometry followed by the radial change \(\rho\mapsto\rho/\sqrt{2-\rho^2}\), where \(\rho=|a|\). Each open unit disk is in turn radially diffeomorphic to its whole vector space.

Near a collision take the tubular coordinates from \(e_p(w)=\sqrt{1-|w|^2}\,p+w\). Then \(p=a/|a|\), \(w=(x-y)/2=tv\) and \(t^2=1-|a|^2\). On the exceptional divisor the limiting map is
\[
 (p,\ell)\longmapsto[v+iJ_pv],\qquad v\in\ell,\ |v|=1.
\]
Every isotropic line admits a representative \(v+ib\) with \(v,b\) orthonormal, unique up to phase. Its inverse image is \(p=v\times b\) with \(\ell=\operatorname{span}_{\mathbb R}(v,b)\subset T_pS^6\). Indeed \(p\times v=b\), so \(\ell\) is a \(J_p\)-complex line. Phase rotates the orthonormal pair through an oriented basis of the same plane and preserves \(v\times b\). Conversely \(v\times J_pv=p\). Local smooth choices of unit representative show that both maps are smooth, giving a diffeomorphism \(\mathbb{P}(TS^6,J)\to Q\).

This is the oriented-plane description of the divisor in Bryant \cite[p.~200]{bryant} and Kotschick--Thung \cite[Proposition~8]{kotschick-thung}. For \(0<\rho<1\), the pair determines the flag \(L\subset P\), where \(P\) is the plane oriented by \((v,J_pv)\). Conversely that flag determines \(p=v\times b_0\), where \((v,b_0)\) is an oriented orthonormal basis of \(P\) with \(v\in L\). At \(\rho=0\) the antipodal pair remembers only \(L\in\mathbb{RP}^6\). At \(\rho=1\) the point of the quadric remembers only \(P\). These are the two projections of the double mapping cylinder in \cite[Theorem~1.5]{le-seade-verjovsky}.

To check transverse smoothness, choose a local unit frame \(v_0\) of \(\ell\), write \(v=e^{J_p\theta}v_0\), and put \(A=v_0+iJ_pv_0\). The filling coordinate is \(\sigma=t^2e^{2i\theta}\). If \(r=|a|/\sqrt{2-|a|^2}\), then
\[
 \Phi=\left[A+\gamma\bar A\right],\qquad
 \gamma=e^{2i\theta}\frac{1-r}{1+r}
       =\frac{\sigma}{1+\sqrt{1-|\sigma|^2}} .
\]
Indeed \(v+iJ_pv=e^{-i\theta}A\) and \(r=\sqrt{(1-t^2)/(1+t^2)}\), giving the stated expression after rationalization. Under \(v_0\mapsto e^{J_p\varphi}v_0\), the vector \(A\) acquires the factor \(e^{-i\varphi}\), while \(\sigma\) and \(\gamma\) acquire \(e^{-2i\varphi}\). Thus the projective expression also agrees on overlaps at zero.

The function \(\gamma\) is smooth in the real and imaginary parts of \(\sigma\) and \(d\gamma_0=\tfrac12\,\mathrm{id}\). At a fixed \(\beta=(p,\ell)\) choose a complex linear form \(L\) with \(L(A)\ne0\). A local defining function for \(Q\) is \(f([z])=(z\cdot z)/L(z)^2\). Since \(A\cdot A=0\) and \(A\cdot\bar A=2\), we have
\[
  (f\circ\Phi)(\beta,\sigma)=
    \frac{4\gamma}{L(A+\gamma\bar A)^2},\qquad
  d_\sigma(f\circ\Phi)_{(\beta,0)}=
    \frac{2}{L(A)^2}\,\mathrm{id}.
\]
Thus \(d\Phi\) induces an isomorphism on the real normal two-planes. Its restriction to the tangent space of the zero section is an isomorphism onto \(TQ\) by the preceding diffeomorphism. Hence \(\Phi\) is a local diffeomorphism across the entire divisor. It is bijective there and on the complementary open set, with disjoint images \(Q\) and \(\mathbb{CP}^6\setminus Q\). Its local inverses therefore give a global smooth inverse. Combining this with the smooth family of fillings and the Chern pairing in Section~\ref{sec:hilbert} proves Theorem~\ref{thm:square}.
\subsection{What this gives for the proposed fourfold construction}
Use the complex six-sphere \(M\) from the specific fibration used by Liang, Liu and Zhang \cite{liang}. It carries a holomorphic involution \(\iota\) with fixed locus a curve \(C\simeq\mathbb{CP}^1\) and with both normal eigenvalues equal to \(-1\). Their exotic smooth \(\mathbb{CP}^3\) with \(c_1=0\) is the complex manifold
\[
 Y=(\operatorname{Bl}_C M)/\langle\iota\rangle .
\]
It occurs naturally in \(H=M^{[2]}\). Away from \(C\), send an orbit to the unordered pair \(\{x,\iota x\}\). In \(M^2\) the graph of \(\iota\) meets the diagonal along \(C\). In linearizing coordinates at a fixed point, \(\iota(z,u,v)=(z,-u,-v)\). The diagonal ideal on the graph is \((u,v)\), so its strict transform in \(\operatorname{Bl}_\Delta(M^2)\) is \(\operatorname{Bl}_C M\). It is invariant under exchange, which acts on the graph as \(\iota\). Taking the finite quotient yields a closed holomorphic embedding of \(Y\) in \(H\); invariants preserve the quotient of the local structure sheaf by the invariant ideal by averaging.

\textbf{Obstruction to an intermediate fourfold in this sixfold.} In the fibration \(f:M\to\mathbb{CP}^1\) of \cite{alpoge}, a very general smooth fiber is a complex two-torus whose N\'eron--Severi group is generated by a class of negative square \cite[Corollary 9.8 and Proposition 9.5]{alpoge}. Such a torus contains no compact curve: an effective curve class is nonzero by pairing with a K\"ahler form, and has nonnegative square by intersecting a general translate. Every compact curve in \(M\) is vertical, since \(c_1(f^*\mathcal{O}(1))=0\) in \(H^2(M;\mathbb{Z})=0\). Thus through a very general point of \(M\) there is no compact curve.

Suppose a compact irreducible analytic fourfold \(W\subset H\) contained \(Y\). The inverse image of \(W\) under the finite quotient is the union of a component dominating \(W\) and its translate under exchange, since exchange acts transitively on each quotient fiber. The irreducible strict transform of the graph lies in one of these components. Blowing that component down gives a compact irreducible analytic set \(V\subset M^2\) containing the graph of \(\iota\). It has dimension four because blow-down is an isomorphism off the exceptional divisor and the graph meets this open set.

The first projection \(V\to M\) is proper and surjective. By the fiber dimension theorem, outside a proper analytic subset of \(M\) its fibers have dimension \(4-3=1\). At the particular point \((x,\iota x)\), the fiber is cut out in the pure four-dimensional space \(V\) by three local coordinate functions from \(M\). Every component of the fiber germ there consequently has dimension at least one. For general \(x\) an irreducible component through \((x,\iota x)\) is therefore a compact curve in \(\{x\}\times M\). We can choose such an \(x\) with \(\iota x\) in a very general, curve-free fiber of \(f\), a contradiction. This excludes an intermediate fourfold in the specified embedding. A different embedding of \(Y\) in a fourfold remains possible.

\subsection{Recovering the complex threefold}\label{sec:recover}
The preceding construction embeds a particular complex threefold in \(H\). A complementary question is whether we can recover the base manifold from its Hilbert square, i.e. is it possible that non-isomorphic base manifolds have isomorphic Hilbert squares. Under a numerical condition on the canonical bundles, we provide the folowing answer to this question:

\begin{theorem}\label{thm:recover}
Let \(X\) and \(Y\) be connected compact complex manifolds of the same dimension \(d\ge3\). Suppose that \(K_X\) and \(K_Y\) have nonnegative degree on every compact irreducible curve in their respective manifolds. Then
\[
 X^{[2]}\cong_{\rm bihol}Y^{[2]}
 \quad\Longleftrightarrow\quad
 X\cong_{\rm bihol}Y.
\]
\end{theorem}

\begin{proof}
Write \(H_X=X^{[2]}\), \(B_X=\operatorname{Bl}_{\Delta}(X\times X)\), and let \(p:B_X\to X\times X\) and \(q:B_X\to H_X\) be the blow-up and double quotient. Denote the exceptional divisors by \(\widetilde E_X\subset B_X\) and \(E_X\subset H_X\). The blow-up canonical bundle formula and the simple branching of \(q\) along \(\widetilde E_X\) give
\[
 K_{B_X}=p^*K_{X\times X}+(d-1)\widetilde E_X
       =q^*K_{H_X}+\widetilde E_X,
 \qquad
 q^*K_{H_X}=p^*K_{X\times X}+(d-2)\widetilde E_X.
\]
Let \(C\subset H_X\) be a compact irreducible curve not contained in \(E_X\), and choose an irreducible component \(C'\) of \(q^{-1}(C)\) that maps onto \(C\). The projections of \(p(C')\) to \(X\) are either points or compact curves, so the hypothesis on \(K_X\) gives \(p^*K_{X\times X}\cdot C'\ge0\). Also \(\widetilde E_X\cdot C'\ge0\), since \(\widetilde E_X\) is effective and does not contain \(C'\). Consequently \(K_{H_X}\cdot C\ge0\). Notice that \(C'\) need not be rational even when \(C\) is rational; this is why the numerical hypothesis concerns all compact curves.

For a line \(\ell\) in a fiber of \(E_X=\mathbb P(TX)\to X\), the normal bundle is \(\mathcal O_{E_X}(-2)\), and adjunction gives \(K_{H_X}\cdot\ell=2-d<0\). Every point of \(E_X\) lies on such a line. Thus \(E_X\) is exactly the union of the rational curves on which \(K_{H_X}\) has negative degree. Any biholomorphism \(\Phi:H_X\to H_Y\) therefore carries \(E_X\) onto \(E_Y\).

Let \(F\cong\mathbb{CP}^{d-1}\) be a fiber of \(\pi_X:E_X\to X\), and set \(f=(\pi_Y\circ\Phi)|_F:F\to Y\). Since \(\pi_Y\circ\Phi|_{E_X}\) is a holomorphic submersion, its differential makes \(f^*TY\) a quotient of \(TE_X|_F\). The tangent sequence of \(\pi_X\) splits over \(F\), because \(H^1(F,TF)=0\), so \(TE_X|_F\cong TF\oplus\mathcal O_F^d\). Hence \(f^*TY\) is globally generated. Its determinant has nonnegative degree on every line in \(F\); the hypothesis on \(K_Y\) gives the opposite inequality. Its first Chern class is therefore zero. A globally generated bundle on projective space with trivial determinant is trivial: a basis of its fiber at one point extends to global sections whose determinant is a nonzero constant. Thus \(f^*TY\) is trivial. Since \(H^0(F,\Omega_F^1)=0\), the differential \(df:TF\to f^*TY\) vanishes, and \(f\) is constant.

The map \(\pi_Y\circ\Phi|_{E_X}\) consequently descends to a holomorphic map \(g:X\to Y\); local holomorphic sections of \(\pi_X\) show holomorphicity. Applying the same argument to \(\Phi^{-1}\) gives a holomorphic inverse. The converse follows from functoriality of the Hilbert square.
\end{proof}

For a complex threefold diffeomorphic to \(S^6\), \(c_1(K_M)=0\) because \(H^2(S^6;\mathbb Z)=0\). Thus the theorem applies without a K\"ahler assumption. In particular, any two nonbiholomorphic complex six-spheres give nonbiholomorphic complex structures on smooth \(\mathbb{CP}^6\), even when their Hodge numbers agree.

\section{Higher symmetric products}\label{sec:symmetric}
\subsection{A local test for higher point spaces}
At a configuration consisting of a double point and \(k-2\) distinct simple points on a complex \(d\)-fold, the Hilbert space locally splits into the length-two model and independent coordinates for the other points. On the exceptional fiber \(F=\mathbb{CP}^{d-1}\), the tangent bundle has contributions \(TF\), a trivial bundle of base directions, and \(\mathcal{O}_F(-2)\). Writing \(t=c_1(\mathcal{O}_F(1))\), this gives \(p_1|_F=(d+4)t^2\).

For \(d\ge3\), if the smooth ambient manifold were \(\mathbb{CP}^{dk}\) and a generator restricted to \(bt\), its known first Pontryagin class would give
\[
 (dk+1)b^2=d+4.
\]
For \(d\ge3,\ k\ge2\), the only solution is \(d=3,\ k=2,\ |b|=1\). Already for three points on a complex six-sphere the equation would read \(10b^2=7\). This argument uses the ordinary neighborhood of a simple collision. Theorem~\ref{thm:degree} also covers resolutions that change that neighborhood.

\subsection{A degree obstruction for all higher symmetric products}\label{higher:degree}
Nakaoka \cite{nakaoka} computed the mod-\(p\) cohomology of symmetric products of spheres. For the degree question below, it suffices to use the sum of the orientation classes of the factors and its first Steenrod reduced power. Both are available for permutation quotients of any closed oriented manifold of even dimension.

\begin{theorem}\label{thm:degree}
Let \(M\) be a closed connected oriented smooth manifold of real dimension \(2d\), \(d\ge2\), and let a subgroup \(G\subseteq S_k\) permute the factors of \(M^k\), \(k\ge2\). A continuous map
\[
       \mathbb{CP}^{dk}\longrightarrow M^k/G
\]
of degree \(1\) or \(-1\) is impossible unless \(d=3\) and \(k=2\). In that exceptional case the quotient must use both permutations. The complex six-sphere construction realizes the exception.
\end{theorem}

\begin{proof}
Choose an odd prime \(p>k\). The quotient \(Z=M^k/G\) is an oriented pseudomanifold: the permutations preserve orientation and act freely on configurations of distinct points. Its top class is the one determined on that open set. Pullback from \(Z\) identifies its cohomology over \(\mathbb F_p\) with the \(G\)-invariants on \(M^k\). Indeed, an equivariant triangulation can be subdivided so that simplex stabilizers fix their simplices pointwise; quotient cochains are invariant cochains, and averaging is exact because \(p\nmid |G|\).

Let \(u_i\) be the pullback of the mod-\(p\) orientation class of \(M\) from factor \(i\). Their sum descends to a class \(\alpha\in H^{2d}(Z;\mathbb F_p)\). Its \(k\)-th power pulls back to \(k!\,u_1\cdots u_k\ne0\). A degree-\(\pm1\) map therefore pulls \(\alpha\) back to \(c x^d\), where \(x\) is the hyperplane generator modulo \(p\) and \(c\ne0\).

Apply the first Steenrod reduced power \(P^1\), which raises degree by \(2(p-1)\) \cite{steenrod}. We have \(P^1u_i=0\) by dimension on the individual factor. Additivity, naturality and injectivity of quotient pullback give \(P^1\alpha=0\). On projective space the Cartan formula gives
\[
           P^1(c x^d)=c d\,x^{d+p-1}.
\]
This is nonzero whenever \(p\nmid d\) and \(p\le d(k-1)+1\).

For \(k\ge3\), Bertrand's postulate supplies a prime strictly between \(\max(d,k)\) and \(2\max(d,k)\). It meets both conditions: if \(d\ge k\), then \(2d\le d(k-1)\), while if \(k\ge d\), then \(2k-1\le d(k-1)+1\). For \(k=2,\ d\ge4\), take a prime strictly between \(\lfloor d/2\rfloor\) and \(2\lfloor d/2\rfloor\). It lies strictly between \(d/2\) and \(d\), so it does not divide \(d\). For \(d=k=2\), take \(p=3\). This leaves only \((d,k)=(3,2)\). There the trivial subgroup is impossible: the two factor orientation classes in \(H^6(M^2;\mathbb Q)\) are independent, whereas a nonzero-degree map into \(M^2\) must inject its rational cohomology into that of \(\mathbb{CP}^6\).
\end{proof}

Theorem~\ref{thm:degree} excludes maps of degree \(\pm1\) except in the Hilbert-square case. For a threefold symmetric product, the same reduced-power calculation yields the stronger divisibility statement below for maps of arbitrary degree.

\begin{proposition}\label{prop:degree-three}
Let \(M\) be a closed connected oriented smooth six-manifold, and let \(Z\) be a closed connected oriented smooth eighteen-manifold whose integral cohomology ring is that of \(\mathbb{CP}^9\). For every continuous map \(f:Z\to\operatorname{Sym}^3(M)\),
\[
                 5^3 7^3=42875\ \mid\ \deg f.
\]
The orientation of the symmetric product is determined by the ordered distinct-point configurations.
\end{proposition}

\begin{proof}
Fix \(p=5\) or \(7\) and work with \(\mathbb Z_{(p)}\)-coefficients. The quotient map \(\rho:M^3\to\operatorname{Sym}^3(M)\) has degree \(6\), which is invertible in \(\mathbb Z_{(p)}\). Transfer identifies the cohomology of the quotient with the invariant part of the cohomology of \(M^3\). If \(u_i\) is the pullback of the orientation class from the \(i\)-th factor, there is a class \(\alpha\in H^6(\operatorname{Sym}^3(M);\mathbb Z_{(p)})\) with \(\rho^*\alpha=u_1+u_2+u_3\). Since \((u_1+u_2+u_3)^3=6u_1u_2u_3\), integration and \(\deg\rho=6\) give \(\int_{\operatorname{Sym}^3(M)}\alpha^3=1\).

Choose a generator \(x\in H^2(Z;\mathbb Z)\). We may write \(f^*\alpha=c x^3\), with \(c\in\mathbb Z_{(p)}\), and integration gives \(\deg f=\pm c^3\) in \(\mathbb Z_{(p)}\). The first mod-\(p\) Steenrod reduced power of each \(u_i\) vanishes by dimension. Quotient pullback is injective modulo \(p\), hence \(P^1\bar\alpha=0\). Naturality and the Cartan formula now give
\[
 0=P^1(f^*\bar\alpha)=3\bar c\,\bar x^{p+2}.
\]
For \(p=5,7\), the power \(\bar x^{p+2}\) is nonzero in \(H^*(Z;\mathbb F_p)\), and \(3\) is invertible. Thus \(c\in p\mathbb Z_{(p)}\), so \(p^3\) divides \(\deg f\). Applying both primes gives the assertion.
\end{proof}

Composing a resolution of the component containing the distinct configurations with the Hilbert--Chow map gives a proper generically one-to-one map to \(\operatorname{Sym}^k M\). Any resolution of the symmetric product itself also has degree one. None can have the smooth type of \(\mathbb{CP}^{dk}\) when \(dk>6,\ d,k\ge2\). For curves, \(\operatorname{Sym}^k\mathbb{CP}^1=\mathbb{CP}^k\) in every length, as seen from binary forms of degree \(k\). A positive-genus curve contributes nonzero \(H^1\) to its symmetric product.

\subsection{What the higher Hilbert spaces actually produce}
For a complex threefold \(M\) diffeomorphic to \(S^6\), let \(M^{[k]}\) denote the full length-\(k\) point space. Its Euler characteristics begin as follows.
\[
\begin{array}{c|rrrr}
 k&1&2&3&4\\ \hline
 \dim_{\mathbb C}\text{ of the distinct-point component}&3&6&9&12\\
 \chi_{\rm top}(M^{[k]})&2&7&18&47\\
 \chi_{\rm top}(\mathbb{CP}^{3k})&4&7&10&13
\end{array}
\]
To see this without any projectivity hypothesis on \(M\), stratify by the distinct support points and their multiplicities. A support of multiplicity \(j\) contributes the Euler characteristic \(p_3(j)\) of the punctual length-\(j\) Hilbert scheme of \(\mathbb C^3\). This projective punctual scheme has one torus-fixed point per monomial ideal, so \(p_3(j)\) counts plane partitions of \(j\). Local analytic coordinates identify the punctual fibers on each ordered support stratum. Compactly supported Euler characteristic is multiplicative on these locally trivial families, and a finite free quotient divides it by the covering degree.

The Euler characteristic of \(s\) ordered distinct points is the falling factorial \(\chi(M)(\chi(M)-1)\cdots(\chi(M)-s+1)\), by successively forgetting the last point. Since \(\chi(M)=2\), only one or two supports contribute. The resulting series, using MacMahon's formula, is
\[
 \sum_{k\ge0}\chi_{\rm top}(M^{[k]})q^k
 =\left(\sum_{j\ge0}p_3(j)q^j\right)^2
 =\prod_{j\ge1}(1-q^j)^{-2j}.
\]
Here \(p_3(0)=1\). For \(k\ge3\) the first two factors alone give a lower bound \((k+1)+4(k-1)>3k+1\). The support argument explains the analytic case; \cite{gusein} gives the corresponding general generating-series construction for smooth quasi-projective varieties. Theorem~\ref{thm:degree} treats resolutions separately.

\textbf{Squaring a second time.} For a compact complex \(d\)-fold \(X\) with Euler characteristic \(e\), distinct pairs contribute \(e(e-1)/2\) and doubled points contribute \(de\). Hence \(\chi(X^{[2]})=e(e+2d-1)/2\). Taking \(X=H\), we get a smooth complex twelvefold of Euler characteristic \(63\). Smooth \(\mathbb{CP}^{12}\) has Euler characteristic \(13\).

\textbf{Scope.} These obstructions use the existence of a degree-one map to the permutation quotient. They leave open complex structures on higher projective spaces constructed in other ways. Sections~\ref{sec:neighborhoods} and~\ref{sec:extensions} keep the exotic sixfold itself as a complex submanifold and ask which ambient structures its tangent bundle permits.

\section{Complex neighborhoods of the sixfold}\label{sec:neighborhoods}
\subsection{Chern classes, symplectic taming, and hypersurfaces}\label{higher:seven}
Recall \(h=c_1(TH)\), fixed by its pairing with the exceptional line in Section~\ref{sec:hilbert}. On \(E=\mathbb{P}(TM)\), put \(\xi=c_1(\mathcal{O}_E(1))\). Then \(h|_E=\xi\), \([E]=-2h\), and \(\xi^3=-\pi^*c_3(TM)\) by the projective bundle formula. In particular \(\int_E\xi^5=-2\) and \(\int_Hh^6=1\).

The relative Euler sequence of \(E\to M\) and its normal line give
\[
 c(TH)|_E=(1-\xi^3)(1+3\xi+3\xi^2)(1-2\xi).
\]
Indeed \(c(\pi^*TM\otimes\mathcal{O}_E(1))=(1+\xi)^3-\xi^3\). Restriction to \(E\) detects powers of \(h\) through degree ten. The top Chern number is \(\chi(H)=7\). These facts give the coefficients of \(c(TH)\) in ascending powers of \(h\):
\[
 (1,1,-3,-7,-1,3,7).
\]
Equivalently, for \(A=K_H^{-1}\),
\[
 [TH]=6A+A^{-1}+A^{-2}-A^2-1\quad\text{in }K^0(H).
\]
The virtual bundle on the right has Chern series \((1+h)^6(1-h)(1-2h)/(1+2h)\) with exactly these coefficients through degree six. Its rank is six, and \(K^0(H)\simeq K^0(\mathbb{CP}^6)\) is torsion-free, so equality of Chern characters suffices.

\begin{theorem}\label{thm:not-tamed}
The complex structure on \(H=M^{[2]}\) is not tamed by any symplectic form.
\end{theorem}

\begin{proof}
Suppose that a symplectic form \(\omega\) tames the complex structure \(J\) of \(H\). Since \(H\) is diffeomorphic to \(\mathbb{CP}^6\) and \(h=c_1(TH)\) is primitive, \([\omega]=a h\) for some real \(a\). The complex line \(\ell\) in a fiber of \(E\to M\) has \(\int_\ell h=1\), and taming gives \(\int_\ell\omega>0\); hence \(a>0\). The restriction \(\omega|_E\) tames the complex structure on the compact complex fivefold \(E\), so its fifth power is positive on the complex orientation and \(\int_E\omega^5>0\). But the divisor class \([E]=-2h\) and \(\int_Hh^6=1\) yield
\[
 \int_E\omega^5=\int_H[E]\wedge[\omega]^5=-2a^5<0,
\]
a contradiction.
\end{proof}

\begin{proposition}\label{prop:hypersurface}
There is no embedding of \(H\) as a complex hypersurface in any complex structure on smooth \(\mathbb{CP}^7\). The same obstruction holds for an almost complex ambient structure that induces the given complex structure on \(H\).
\end{proposition}

\begin{proof}

Suppose \(j:H\hookrightarrow X\) were such an embedding. Write \(j^*x=b h\) for an integral generator \(x\) and \(\epsilon=\int_Xx^7=\pm1\), using the complex orientation. Poincar\'e duality gives the divisor class \(\epsilon b^6x\) and normal Chern class \(\epsilon b^7h\). Restricting \(p_1\) gives \(8b^2=7+b^{14}\), whence \(|b|=1\). Choose \(x\) so that \(b=1\), and let \(B\) be a topological complex line of first Chern class \(x\). Then \(B|_H\simeq A\) smoothly and \(N_{H/X}\simeq B^\epsilon|_H\).

Consider \(W_\epsilon=6B+B^{-1}+B^{-2}-B^2-1+B^\epsilon\). By the tangent-normal sequence, its restriction agrees with \(TX|_H\) in \(K^0(H)\). In \(K^0(\mathbb{CP}^7)=\mathbb{Z}[B-1]/((B-1)^8)\), the kernel of restriction to \(H\) is generated by \((B-1)^7\). Its degree-fourteen Chern character is \(x^7\). Newton's identity, with the lower Chern classes fixed, thus forces
\[
       \int_X c_7(TX)\equiv\int_X c_7(W_\epsilon)\pmod{6!}.
\]
The Chern series of \(W_\epsilon\) is \((1+x)^6(1-x)(1-2x)(1+\epsilon x)/(1+2x)\). Its \(x^7\) coefficient is \(2\) for \(\epsilon=1\) and \(-12\) for \(\epsilon=-1\). The two candidate Chern \emph{numbers} are therefore \(2\) and \(12\). Neither is congruent to \(\chi(\mathbb{CP}^7)=8\) modulo \(720\).
\end{proof}

\subsection{Integrable neighborhoods and compact complex blow-ups}\label{higher:open}
Fix \(n\ge7\), put \(r=n-6\), and keep \(A=K_H^{-1}\), with \(c_1(A)=h\). The complex manifold
\[
                    U_n=\operatorname{Tot}(A^{\oplus r})\longrightarrow H
\]
is diffeomorphic to \(\mathbb{CP}^n\setminus\mathbb{CP}^{\,r-1}\). To check the smooth type, identify \(H\) with \(\mathbb{CP}^6\) so that \(h\) corresponds to the hyperplane generator, using conjugation if necessary. Smooth complex lines are classified by \(c_1\). In standard projective space, projection away from \(\mathbb{P}(\mathbb C^r)\) sends \([z:w]\) to \([z]\in\mathbb{CP}^6\); \(w\) represents a linear map \(\mathbb Cz\to\mathbb C^r\). This identifies the complement with the total space of \(\mathcal{O}(1)^{\oplus r}\).

The bundle retracts onto \(H\). Its tangent sequence gives \(c_1(TU_n)=(r+1)h=(n-5)h\). On the standard open manifold the first Chern class is \((n+1)h\). Thus in every dimension \(n>6\) this gives a nonstandard \emph{integrable} structure on the open complement. With the displayed identification, a smooth integrable extension agreeing on that entire complement is impossible: fixing one line in \(\mathbb C^r\) gives a copy of \(\operatorname{Tot}(A)\) whose closure is a smooth linear \(\mathbb{CP}^7\). Its tangent planes are complex on a dense set and would stay complex everywhere. It would contain \(H\) as a complex hypersurface, contrary to Proposition~\ref{prop:hypersurface}.

\textbf{A compact result in every dimension \(n>6\).} Projectivizing with the convention of lines gives
\[
 Z_n=\mathbb{P}_H(\mathcal{O}\oplus A^{\oplus r})
 \ \cong_{\rm diff}\ \operatorname{Bl}_{\mathbb{CP}^{r-1}}\mathbb{CP}^n .
\]
For the standard base, the graph closure of the same linear projection is the blow-up of its center and is \(\mathbb{P}_{\mathbb{CP}^6}(\mathcal{O}(-1)\oplus\mathcal{O}^{\oplus r})\). Tensoring by \(\mathcal{O}(1)\) gives the displayed model. Projective bundle fibers have no higher cohomology of \(\mathcal{O}\), so \(h^{0,1}(Z_n)=h^{0,1}(H)\). Viaclovsky's family \cite{viaclovsky} therefore gives infinitely many pairwise nonbiholomorphic complex structures on each of these fixed compact smooth blow-ups. Here \(b_2(Z_n)=2\).

The divisor at infinity is \(D=H\times\mathbb{CP}^{r-1}\). Even the expected holomorphic blow-down fails. Suppose a map to a smooth \(n\)-fold contracted each \(H\)-factor of \(D\) to a point and were biholomorphic away from \(D\). Its Jacobian would have divisor \(mD\) with \(m>0\). On \(H\times\{q\}\), the normal of \(D\) is \(A^{-1}\), while the other normal directions in \(D\) are trivial. Adjunction gives \(c_1(K_{Z_n}|_H)=c_1(K_H\otimes A)=0\). The Jacobian formula on this same contracted slice gives \(-mh\), a contradiction.

\subsection{An obstruction to compactification with analytic boundary}\label{higher:cap}
\begin{theorem}\label{thm:compactification}
For \(n=6+r\), \(r\ge1\), the open complex manifold \(U_n=\operatorname{Tot}(A^{\oplus r})\) cannot occur as \(X\setminus D\), where \(X\) is a compact complex manifold diffeomorphic to \(\mathbb{CP}^n\) and \(D\) is a closed complex analytic subset. This allows an arbitrary identification of the open set and an arbitrary analytic boundary.
\end{theorem}

\begin{proof}

\textbf{The topology forces the boundary to have complex codimension at least seven.}
The inclusion \(U_n\hookrightarrow X\) pulls \(c_1(TX)\) back to \((r+1)h\ne0\). It therefore induces an isomorphism in degree two over \(\mathbb Q\). Powers of a degree-two generator give all the cohomology of \(U_n\simeq H\simeq\mathbb{CP}^6\), so restriction is an isomorphism through degree twelve. Both spaces have vanishing \(H^{13}\). The exact sequence of the pair gives \(H^j(X,U_n;\mathbb Q)=0\) for \(j\le13\). By Poincar\'e--Alexander duality this is \(H_{2n-j}(D;\mathbb Q)\). A compact analytic set of complex dimension \(s\) has a nonzero fundamental class in degree \(2s\), obtained from any component of maximal dimension. Hence \(\dim_{\mathbb C}D\le n-7\).

\textbf{Cohomology of the open set.} Scheja's extension theorem \cite{scheja} for the structure sheaf on a smooth complex manifold identifies \(H^q(X,\mathcal{O}_X)\) with \(H^q(U_n,\mathcal{O}_{U_n})\) for \(0\le q\le5\) in this codimension. These groups are finite-dimensional by compactness of \(X\).
Let \(\pi:U_n\to H\) be the bundle projection. Locally its preimage over a trivializing coordinate polydisc is Stein, so \(H^q(U_n,\mathcal{O}_{U_n})=H^q(H,\pi_*\mathcal{O}_{U_n})\) by Cartan's theorem B and Leray. For each \(N\), the sheaf \(\bigoplus_{k=0}^N A^{-k}\) is a direct summand of \(\pi_*\mathcal{O}_{U_n}\) as an \(\mathcal{O}_H\)-module. Include the fiber polynomials in the first copy of \(A\). To split this inclusion, set the other coordinates to zero and take the Taylor coefficients through degree \(N\). These operations commute with changes of line frame. The split injections on cohomology, for every \(N\), force
\[
             H^q(H,A^{-k})=0\qquad (0\le q\le5,\ k\gg0).
\]

\textbf{The contradiction on the sixfold.} Since \(K_H=A^{-1}\), Serre duality and Hirzebruch--Riemann--Roch give, for large \(k\),
\[
 h^0(H,A^{k-1})=h^6(H,A^{-k})=\chi(H,A^{-k})
             =\frac{k^6}{6!}\int_H h^6+O(k^5)=\frac{k^6}{720}+O(k^5).
\]
Thus \(A\) is big and \(H\) is Moishezon \cite{ma-marinescu}. A smooth compact Moishezon manifold admits a smooth projective modification and satisfies \(b_1=2h^{0,1}\), by bimeromorphic invariance of both numbers. Since \(b_1(H)=0\), this would give \(h^{0,1}(H)=0\).

On a complex six-sphere, Gray \cite[Theorem~5]{gray} proved \(h^{0,1}(M)>0\). More precisely, \(h^{0,1}(M)=1+h^{0,2}(M)\), as recorded by Angella \cite[eq.~(3)]{angella}. To see this, a holomorphic three-form \(\eta\) is closed and hence exact because \(H^3(M;\mathbb C)=0\). Stokes' theorem gives \(\int_M i\,\eta\wedge\bar\eta=0\), so \(\eta=0\). Serre duality gives \(h^{0,3}(M)=0\), and Riemann--Roch gives \(\chi(M,\mathcal O_M)=\int_M c_1c_2/24=0\). Expanding this Euler characteristic gives the stated identity. The equality \(h^{0,1}(H)=h^{0,1}(M)\) from Section~\ref{sec:hilbert} now yields the contradiction.
\end{proof}

\section{Almost complex extensions}\label{sec:extensions}
\subsection{Even dimensions}\label{higher:even}
\begin{theorem}\label{thm:even}
For every even \(n\ge8\) and every open neighborhood \(V\) of the deleted linear \(\mathbb{CP}^{n-7}\), there is a smooth almost complex structure \(J_V\) on \(\mathbb{CP}^n\) such that
\[
 c_1(T\mathbb{CP}^n,J_V)=(n-5)x,\qquad
 J_V|_{\mathbb{CP}^n\setminus V}=J_{U_n}|_{\mathbb{CP}^n\setminus V}.
\]
Here \(x\) is the standard integral generator identified with \(h\) on \(H\). The choice of \(J_V\) may depend on \(V\).
\end{theorem}

\begin{proof}

Almost complex structures on projective spaces were studied by A. Thomas \cite{thomas}. Here the additional requirement is to match the integrable structure near \(H\). Let \(B\) be the line of first Chern class \(x\). Define a stable class of complex rank \(n\) by
\[
 W_n=nB+B^{-1}+B^{-2}-B^2-1.
\]
Its realification agrees stably with \(T\mathbb{CP}^n\): its difference from the standard complex tangent class \((n+1)B-1\) is a sum of differences of lines and their conjugates. Its restriction to \(H\) is \([TH]+rA\) by Section~\ref{higher:seven}. For a formal variable \(t\) representing \(x\), its total Chern series is
\[
 (1+t)^n\frac{(1-t)(1-2t)}{1+2t}
      =(1+t)^n\left(t-2+\frac{3}{1+2t}\right).
\]
Its linear coefficient is \(n-5\), and its top coefficient is \(n-2+3(-1)^n\) by the binomial theorem. For even \(n\) this equals \(n+1\), the Euler number. The stable-to-unstable criterion of E. Thomas now gives an actual almost complex structure realizing \(W_n\) on the oriented tangent bundle \cite[Theorem 2.2]{mills}. For instance, on \(\mathbb{CP}^8\) it has \(c_1=3x\) and \(c_8=9x^8\).

\textbf{Matching the two structures along \(H\).} Denote the global structure by \(J_0\), the integrable structure on \(U_n\) by \(J_1\), and the underlying real tangent bundle restricted to \(H\) by \(T\). Their complex orientations on \(T\) agree. Indeed, \(\int_Hh^6=1\) gives the standard orientation on the base under the identification in Section~\ref{higher:open}, and the chosen identification of the normal bundles is complex linear. The two complex bundles \((T,J_i)\) represent the same class in \(K^0(H)\). They have rank \(n\ge8\) over a real \(12\)-manifold, so complex cancellation gives a bundle isomorphism \(g:(T,J_0)\to(T,J_1)\). Here \(BU(n)\to BU\) is \((2n+1)\)-connected by the sphere fibrations for the unitary groups.

As a real bundle map, \(g\) is an orientation-preserving automorphism of \(T\). Choose a real complement \(F\) so that \(T\oplus F\) is trivial. In a trivialization, \(g\oplus 1_F\) defines a map to \(GL^+(N,\mathbb R)\), where \(N=\operatorname{rank}(T\oplus F)\). Its stable class lies in
\[
 [H,SO]=\widetilde{KO}^{-1}(\mathbb{CP}^6)=0
\]
by Fujii's computation \cite[proof of Proposition~2.1]{mills}. Enlarging \(F\) by a trivial summand if necessary, we obtain a homotopy of \(g\oplus 1_F\) to the identity. The block inclusion of oriented automorphism bundles for \(T\) and \(T\oplus F\) is \((2n-1)\)-connected on each fiber. This follows from polar decomposition and \(SO(j+1)/SO(j)=S^j\) for \(j\ge2n\). Relative obstruction theory over \((H\times[0,1],H\times\{0,1\})\) therefore deforms this homotopy, fixing its ends, to one in the smaller automorphism bundle, since \(13=\dim(H\times[0,1])<2n-1\). Conjugating \(J_0\) by the resulting path from the identity to \(g\) gives a homotopy to \(J_1\) on \(H\). It can be smoothed relative to its ends.

\textbf{Agreement outside \(V\).} Choose a smooth connection on the ambient real tangent bundle over \(U_n\). For \(v\in U_n\), evaluate each \(J_i\) at \(s v\) and transport it back to \(v\) along the fiber segment, \(0\le s\le1\). This deforms \(J_i\) to the radial transport of its restriction to the zero section. Inserting the transported homotopy on \(H\) between these two paths gives a smooth path \(K_t\) from \(J_0\) to \(J_1\) on \(U_n\), after making the parameters constant near the joins. Take a smooth function \(\chi:U_n\to[0,1]\) of compact support equal to one near the compact set \(\mathbb{CP}^n\setminus V\). Set \(J_V(v)=K_{\chi(v)}(v)\) on \(U_n\) and use \(J_0\) off its support. Replacing \(\chi\) by \(t\chi\) gives a global homotopy to \(J_0\), so \(c_1(J_V)=(n-5)x\).
\end{proof}

The construction localizes any failure of integrability to \(V\). It gives neither a single structure agreeing on all of \(U_n\) nor a globally integrable structure with this first Chern class.

\subsection{Why the same neighborhood fails in odd dimensions}\label{higher:odd}
\begin{lemma}\label{lem:odd}
For odd \(n\ge5\), if \(J\) is an almost complex structure on smooth \(\mathbb{CP}^n\) and \(c_1(T\mathbb{CP}^n,J)=ax\), then \(a\equiv n+1\pmod4\). Changing the sign of the generator preserves this congruence.
\end{lemma}

\begin{proof}
Write \(y=B-1\) in \(K^0(\mathbb{CP}^n)=\mathbb{Z}[y]/(y^{n+1})\), and let a bar denote complex conjugation. The difference \(\delta=[(T\mathbb{CP}^n,J)]-[(n+1)B-1]\) has zero realification. Complexifying gives \(\delta+\bar\delta=0\). Since \(\bar y=-y/(1+y)\), the first nonzero coefficient \(m y^j\) of an anti-invariant polynomial must have \(j\) odd. If \(j<n\) and its next coefficient is \(b\), comparison in degree \(j+1\) gives \(jm+2b=0\). Thus \(m\) is even. Subtract \((m/2)(y^j-\bar y^j)\) and repeat. This expresses \(\delta\) as an integral sum of differences \(B^i-B^{-i}\) and a multiple of \(y^n\).

The latter term changes only the top Chern coefficient, by a multiple of \((n-1)!\), hence by \(0\) modulo \(4\). For each line difference one has
\[
 c(B^i-B^{-i})=\frac{1+ix}{1-ix}
 \equiv 1+2(i\bmod2)\frac{x}{1-x}\pmod4 .
\]
Products of the corrections vanish modulo \(4\). Comparing first Chern classes therefore gives
\[
 c(T\mathbb{CP}^n,J)\equiv(1+x)^{n+1}
          \left(1+(a-n-1)\frac{x}{1-x}\right)\pmod4 .
\]
The coefficient of \(x^{n-1}\) in \((1+x)^{n+1}/(1-x)\) is \(2^{n+1}-n-2\), an odd integer. With either complex orientation the top Chern coefficient is \(\pm(n+1)\equiv n+1\pmod4\). Comparing top coefficients proves the assertion.
\end{proof}

For odd \(n\ge7\), the integrable neighborhood in Section~\ref{higher:open} has first Chern coefficient \(n-5\), differing from \(n+1\) by six. It cannot agree near its zero section with any almost complex structure on the closed \(\mathbb{CP}^n\). This applies to the specified normal bundle \(A^{\oplus(n-6)}\). At \(n=7\) the stronger result in Section~\ref{higher:seven} also rules out the other sign of the normal line.

\textbf{The remaining analytic question.} In even dimensions the matching almost complex structures have the required real tangent bundle and an integrable neighborhood. For such a structure to be integrable everywhere its Nijenhuis tensor must vanish everywhere. Shrinking the region where it may be nonzero gives no convergence or removability result. In odd dimensions the congruence rules out this particular first Chern class but leaves other classes possible. The arguments here give no nonstandard integrable structure on any \emph{closed} \(\mathbb{CP}^n\) with \(n>6\), and do not prove uniqueness there.

\end{document}